\documentclass[11pt]{article}
\usepackage{amsfonts}
\usepackage{latexsym, amssymb, amsmath, amscd, amsfonts, epsfig, graphicx, colordvi}
\usepackage{ifpdf}
\usepackage{graphicx}
\usepackage{xcolor}
\usepackage[dvipsnames]{pstricks}
\newtheorem{thm}{Theorem}[section]

\newtheorem{pro}[thm]{Proposition}

\newtheorem{defi}[thm]{Definition}

\newtheorem{exam}[thm]{\bf Example}
\def\pf{\noindent{\it Proof.} }
\def\qed{\nopagebreak\hfill{\rule{4pt}{7pt}}
	\medbreak}

\def\qed{\nopagebreak\hfill{\rule{4pt}{7pt}}
	\medbreak}

\makeatletter
\def\ExtendSymbol#1#2#3#4#5{\ext@arrow 0099{\arrowfill@#1#2#3}{#4}{#5}}
\makeatother

\title{Lattice paths and the Rogers--Ramanujan--Gordon type theorems  with parity considerations}
\author {Diane Y. H.
	Shi}
\date{	\vspace{15pt} School of Mathematics, \\ Tianjin University, Tianjin 300072, P.R. China
	\vskip 0.2 cm
	Email: shiyahui@tju.edu.cn}

\begin{document}
	
\maketitle
%--------------------------------------------------------------------
\noindent {\bf Abstract.} Andrews imposed parity restrictions on the Rogers–Ramanujan–Gordon type partitions, yielding fruitful results. These results were later, advanced by Kur\c{s}ung\"{o}z, Kim, and Yee.  In this paper, we  construct a bijection between lattice paths with three kinds of unitary steps and the Rogers--Ramanujan--Gordon type partitions.  We  give some results involving parity considerations on  lattice paths  by this bijection  which are  Andrews' partition results counterpart. We also give some new results on lattice paths.
	
\noindent {\bf Keywords:} lattice path, the Rogers--Ramanujan--Gordon Theorem,
the Andrews--Gordon identity, partition, Gordon marking 

\noindent {\bf AMS Subject Classification:} 05A17, 11P84	

\section{ Introduction}
	
In 2010, Andrews [6] investigated a variety of parity questions in classical partition
identities. In particular, he revisited the Rogers–Ramanujan–Gordon theorem along with
the generating function  and extended his results by considering some additional
restrictions involving parities. In this paper, we will revisit these identities by employing  lattice paths. We firstly recall some related theorems. 

Gordon \cite{gor61} found the following combinatorial generalization of the
Rogers--Ramanujan identities, which has been called the Rogers--Ramanujan--Gordon theorem,
see Andrews \cite{and66}.

\begin{thm}\label{GRR}
	For $k \geq a\geq 1$, let $B_{k,a}(n)$ denote the number of partitions
	of $n$ of the form $\lambda_1+\lambda_2+\cdots+\lambda_s$, where $\lambda_j\geq \lambda_{j+1}$, $\lambda_j-\lambda_{j+2}\geq 2$ and at most $a-1$ of the
	$\lambda_j$ are equal to $1$. Let $A_{k,a}(n)$ denote the number of partitions of $n$ into parts  not congruent to
	$0,\pm a$ modulo $2k+1$. Then for all $n\geq 0$, we have
	\[A_{k,a}(n)=B_{k,a}(n).\]
\end{thm}
Let us give an overview of some definitions of partitions. An unrestricted partition $\lambda$ of a positive integer $n$ is a non-increasing sequence of positive integers $\lambda_1\geq \cdots\geq \lambda_s>0$ such that $n=|\lambda|=\lambda_1+\cdots+\lambda_s$. Let $p(n)$ be the number of partitions of $n$ and as $0$  has only one partition $\emptyset$,  let $p(0)=1$.
Given a partition $\lambda$, let $f_l(\lambda)$  denote the number of occurrences of
$l$ in $\lambda$.

We adopt the common notations in $q$-series as used in Andrews \cite{and76}. Let
\[(a)_\infty=(a;q)_\infty=\prod_{i=0}^{\infty}(1-aq^i),\]
and
\[(a_1,\ldots,a_k;q)_\infty=(a_1;q)_\infty\cdots(a_k;q)_\infty.\]
We also write
\[(a)_n=(a;q)_n=(1-a)(1-aq)\cdots(1-aq^{n-1}).\]

In 1974, Andrews \cite{and74} derived an identity which can be considered as  the
generating function counterpart of the Rogers--Ramanujan--Gordon theorem. It has been called the Andrews--Gordon identity, and it is an analytic generalization of the Rogers--Ramanujan identities with odd moduli.

\begin{thm}For $k\geq a\geq 1$, we have
	\begin{equation}\label{abRRG}
		\sum_{N_1\geq N_2\geq \cdots \geq N_{k-1}\geq 0}\frac{q^{N_1^2+N_2^2+\cdots+N_{k-1}^2+N_a+\cdots+N_{k-1}}}
		{(q)_{N_1-N_2}\ldots(q)_{N_{k-2}-N_{k-1}}(q)_{N_{k-1}}}=
		\frac{(q^a,q^{2k+1-a},q^{2k+1};q^{2k+1})_{\infty}}
		{(q)_{\infty}}.
	\end{equation}
\end{thm}

Andrews showed that both sides of \eqref{abRRG} satisfy the
same recurrence relation with the same initial conditions. 

Bressoud \cite{bre79, Bre80} obtained a Rogers--Ramanujan--Gordon type theorem and
the corresponding Andrews--Gordon type identity with even moduli.
%The following Rogers--Ramanujan--Gordon type theorem for even moduli is due to Bressoud.

\begin{thm}\label{Bressoud}For $k\geq a\geq 1$, let $D_{k,a}(n)$ denote the number of partitions of $n$ of the form $\lambda=\lambda_1+\lambda_2+\cdots+\lambda_s$,
	such that
	\begin{itemize} \item[(i)] $f_1(\lambda)\leq a-1$, \item[(ii)] $f_l(\lambda)+f_{l+1}(\lambda)\leq k-1$, \item[(iii)] if the equality in condition (ii)
		is attained at $l$, then $lf_l(\lambda)+(l+1)f_{l+1}(\lambda)\equiv a-1 \;(\rm{mod}\;2)$.
	\end{itemize} Let $B_{k,a}(n)$ denote the number of
	partitions of $n$ whose parts are not congruent to
	$0,\pm a$ modulo $2k$. Then for all $n\geq 0$, we have
	\[C_{k,a}(n)=D_{k,a}(n).\]
\end{thm}
The generating function form of the above theorem can be stated as follows.

\begin{thm}For $k\geq a\geq 1$, we have
	\begin{equation}\label{AB}\sum_{N_1\geq N_2\geq\cdots\geq
			N_{k-1}\geq0}\frac{q^{N_1^2+N_2^2+\cdots+N_{k-1}^2+N_{a+1}+\cdots+N_{k-1}}
		}{(q)_{N_1-N_2}\cdots(q)_{N_{k-2}-N_{k-1}}(q^2;q^2)_{N_{k-1}}}
		=\frac{(q^a,q^{2k-a},q^{2k};q^{2k})_\infty}{(q)_\infty}.
	\end{equation}
\end{thm}

 Andrews and Bressoud \cite{andbre84} showed that left hand side of \eqref{abRRG} could be
rephrased in terms of lattice paths with two kinds of steps. In 1986, Bressoud \cite{bre87} reinterpreted these in terms of different lattice paths with three kinds of
steps and gave some direct bijections between the objects counted by  $B_{k,a}(n)$, $D_{k,a}(n)$ and the lattice paths.  

The lattice paths we  studied in this paper \cite{bre87,cor07} are  paths in the first quadrant, that start on the $y$-axis, end on the $x$-axis,
and use three kinds of unitary steps:

North–East NE: $(x,y) \rightarrow  (x + 1,y + 1)$,

South–East SE: $(x,y) \rightarrow (x+1,y-1)$,

East E: $(x,0)\rightarrow (x+1,0)$.

The height of a vertex $(x,y)$ corresponds to the $y$-coordinate and  the weight of a vertex is
its $x$-coordinate. An East step can only appear at height $0$. The paths are either empty or end with a South–East
step. A peak $(x,y)$ is a vertex preceded by a North–East step and followed  by a South–East step.   The major index of a path is the sum of the weights of its peaks. Let $k$ and $a$ be positive integers with $a\leq k$. We say that a path satisfies the special $(k,a)$-conditions if it starts at height $k-a$ and its height is less than $k$.

Bressoud \cite{bre87} gave the following result by combinatorial construction. 
\begin{thm}
For $k \geq a\geq 1$, let $E_{k,a}(n)$ be the number of paths of major index $n$	which satisfy the special $(k,a)$-conditions. Then we have 
\begin{equation}\label{eqRRG}
	\sum_{n=0}^\infty E_{k,a}(n)q^n=\sum_{N_1\geq N_2\geq \cdots \geq N_{k-1}\geq 0}\frac{q^{N_1^2+N_2^2+\cdots+N_{k-1}^2+N_a+\cdots+N_{k-1}}}
	{(q)_{N_1-N_2}\ldots(q)_{N_{k-2}-N_{k-1}}(q)_{N_{k-1}}}.
	\end{equation}

\end{thm}
By this generating functions of $E_{k,a}(n)$ and  $B_{k,a}(n)$, he also derived that \[E_{k,a}(n)=B_{k,a}(n).\]
In 2009, Kur\c{s}ung\"{o}z \cite{kur09} gave the  combinatorial  interpretation of  the summation on the left-hand side of \eqref{abRRG} that  can be viewed as the generating function for $B_{k,a}(n)$ by using the Gordon marking of partitions.

Bressoud \cite{bre87} also derived the following result which correspond to $D_{k,a}(n)$.
\begin{thm}
	For $k \geq a\geq 1$, let $\widetilde{E}_{k,a}(n)$ be the number of paths of major index $n$ which satisfy  the special $(k,a)$-conditions and every peak of height $k-1$ has weight congruent to $a-1$ modulo $2$. Then we have 
	\begin{equation}\sum_{k=0}^\infty \widetilde{E}_{k,a}(n)q^n =\sum_{N_1\geq N_2\geq\cdots\geq
			N_{k-1}\geq0}\frac{q^{N_1^2+N_2^2+\cdots+N_{k-1}^2+N_{a+1}+\cdots+N_{k-1}}}{(q)_{N_1-N_2}\cdots(q)_{N_{k-2}-N_{k-1}}(q^2;q^2)_{N_{k-1}}}
		\end{equation}

\end{thm}
Hence by \eqref{AB}, 	\[\widetilde{E}_{k,a}(n)=D_{k,a}(n).\]

In 2010, Andrews \cite{And10}  put some  parity restrictions on Rogers-Rmanujan-Gordon partitions, and derived series results. A part of results are about the partitions enumerated by $B_{k,a}(n)$ which have the addition  restriction that   the number of occurrences of even parts is even or the number of occurrences of odd parts is even. Andrews considered these partition numbers according to the parities of $k$ and $a$. He did not give these numbers in all conditions, though, at the end of his paper \cite{And10}, Andrews left the cases for the other values of $k$ and $a$ as open problems.   In \cite{kur09,kur10}, Kur\c{s}ung\"{o}z obtained the infinite sunmmation form generating functions  by employing Gordon marking of partitions. And then,  in \cite{kim13}, Kim and Yee derived the generating functions in other cases to make the results complete. We summarize their results as follows.
\begin{thm}[Andrews \cite{And10}, Kur\c{s}ung\"{o}z \cite{kur10}, Kim and Yee \cite{kim13} ]\label{andwe} Suppose $k\geq a\geq 1$ are integers. Let $W_{k,a}(n)$
	denote the number of those partitions enumerated by $B_{k,a}(n)$ with the added restriction
	that even parts appear an even number of times. Hence we have  if  $k\equiv a \pmod{2}$
	
	\begin{align}\nonumber&\sum_{n\geq 0}W_{k,a}(n)q^n=\frac{(-q;q^2)_\infty(q^a,q^{2k+2-a},q^{2k+2};q^{2k+2})_\infty}{(q^2;q^2)_\infty}\\&=\sum_{N_1\geq N_2\geq\cdots\geq
			N_{k-1}\geq0}\frac{q^{N_1^2+N_2^2+\cdots+N_{k-1}^2+2N_{a}+2N_{a+2}+\cdots+2N_{k-2}}
		}{(q^2;q^2)_{N_1-N_2}\cdots(q^2;q^2)_{N_{k-2}-N_{k-1}}(q^2;q^2)_{N_{k-1}}}
	\end{align}
	For $k$ and $a$ have different parities, there is 
	\begin{align}\sum_{n\geq0}W_{k,a}(n)q^n&=\label{Wkad}\sum_{N_1\geq N_2\geq\cdots\geq
		N_{k-1}\geq0}\frac{q^{N_1^2+N_2^2+\cdots+N_{k-1}^2+2N_{a}+2N_{a+2}+\cdots+2N_{k-1}}
	}{(q^2;q^2)_{N_1-N_2}\cdots(q^2;q^2)_{N_{k-2}-N_{k-1}}(q^2;q^2)_{N_{k-1}}}\\\nonumber&=\frac{(-q^3;q^2)_\infty(q^{a+1},q^{2k+1-a},q^{2k+2};q^{2k+2})_\infty}{(q^2;q^2)_\infty}
	\\ &\label{Wkad2}\qquad+\frac{q(-q^3;q^2)_\infty(q^{a-1},q^{2k+3-a},q^{2k+2};q^{2k+2})_\infty}{(q^2;q^2)_\infty}.
	\end{align}
	\end{thm}
\eqref{Wkad} was derived by Kur\c{s}ung\"{o}z  by using combinatorial way, as well as \eqref{Wkad2} was derived by Kim and Yee \cite{kim13}.
\begin{thm}[Andrews \cite{And10}, Kur\c{s}ung\"{o}z \cite{kur10}, Kim and Yee \cite{kim13}]\label{andoverw}
		Suppose $k\geq a\geq 1$ with $k$ odd and $a$ even. Let $\overline{W}_{k,a}(n)$ denote the
		number of those partitions enumerated by $B_{k,a}(n)$ with added restriction that odd
		parts appear an even number of times. Then
		\begin{align}\nonumber&\sum_{n\geq0}\overline{W}_{k,a}(n)q^n=\frac{(-q^2;q^2)_\infty(q^a,q^{2k+2-a},q^{2k+2};q^{2k+2})_\infty}{(q^2;q^2)_\infty}
			\\&=\sum_{N_1\geq N_2\geq\cdots\geq
				N_{k-1}\geq0}\frac{q^{N_1^2+N_2^2+\cdots+N_{k-1}^2+n_1+n_3+\cdots+n_{a-3}+N_{a-1}
					+N_a+\cdots+N_{k-1}}
			}{(q^2;q^2)_{n_1}\cdots(q^2;q^2)_{n_{k-2}}(q^2;q^2)_{n_{k-1}}}
		\end{align}
		where $n_i=N_i-N_{i+1}$.

	For $k\ge a\ge 1$ with $k$ even and $a$ odd,
\begin{align}
	\sum_{n\geq 0}\overline{W}_{k,a}(n)q^n&=\label{OWakd}\sum_{N_1\geq N_2\geq\cdots\geq
		N_{k-1}\geq0}\frac{q^{N_1^2+N_2^2+\cdots+N_{k-1}^2+n_1+n_3+\cdots+n_{a-2}
			+N_a+\cdots+N_{k-1}}
	}{(q^2;q^2)_{N_1-N_2}\cdots(q^2;q^2)_{N_{k-2}-N_{k-1}}(q^2;q^2)_{N_{k-1}}} \\  &\qquad=\label{OWakd2}\frac{(-q^2;q^2)_\infty(q^{a+1},q^{2k+1-a},q^{2k+2};q^{2k+2})_\infty}{(q^2;q^2)_\infty}.
\end{align}
	\end{thm}
\eqref{OWakd} was derived by 	Kur\c{s}ung\"{o}z \cite{kur10} combinatrially and  Kim and Yee \cite{kim13} in different way, as well as
\eqref{OWakd2} was derived by  Kim and Yee \cite{kim13}.

In this paper, we shall rephrased \eqref{Wkad} and \eqref{OWakd}  by  lattice path. We recall the definition of relative height of a peak which is defined by Bressound \cite{bre87} and retail by Corteel and  Mallet \cite{cor07}. We refer the definition of Corteel and  Mallet \cite{cor07}.

\begin{defi}
	The relative height of a peak $(x,y)$ is the largest integer $h$ for which we can find
two vertices on the path, $(x'
,y-h)$ and $(x'',y-h)$, such that $x'<x<x''$ and such that between
these two vertices there are no peaks of height $> y$ and every peak of height $y$ has abscissa $\geq x$.
	\end{defi}
Now we give the main results of this parper.
\begin{thm}\label{thmwp1}
	Let $k\geq2$, and $k\geq a\geq 1$, let $P_{k,a}(n)$ denote the number of lattice paths $L$ satisfy the special $(k,a)$-conditions whose major index are  $n$  and   the peaks satisfy that the weight and the relative height have the same parity, i.e., all peaks have odd parities.   Then we have 
	\[W_{k,a}(n)=P_{k,a}(n).\]
By the generating function of $W_{k,a}(n)$,  we also can have that
		\begin{equation}\sum_{n\geq 0}P_{k,a}(n)q^n=\sum_{N_1\geq N_2\geq\cdots\geq
			N_{k-1}\geq0}\frac{q^{N_1^2+N_2^2+\cdots+N_{k-1}^2+2N_{a}+2N_{a+2}+\cdots}
		}{(q^2;q^2)_{N_1-N_2}\cdots(q^2;q^2)_{N_{k-2}-N_{k-1}}(q^2;q^2)_{N_{k-1}}}.
	\end{equation}
\end{thm}	
\begin{thm}\label{thmwp2}
	For $k\geq a\geq 1$, let $\overline{P}_{k,a}(n)$ denote the number of lattice paths $L$ satisfying the special  $(k,a)$-conditions whose major indexes are $n$ and the weight of all peaks
	are even. Then we have 
	\[\overline{W}_{k,a}(n)=\overline{P}_{k,a}(n).\]
\end{thm}	
\begin{thm}\label{thmle}
	For $k\geq a\geq 1$, let $Q_{k,a}(n)$ denote the number of lattice paths satisfying the special $(k,a)$-conditions  whose major indexes are $n$ and  all peaks   with odd weight. Then we have the generating functions of $Q_{k,a}(n)$ according to the parities of $k,\ a$  as follows:
\begin{itemize}
	\item If $a$ is odd, then 
\begin{equation}
	\sum_{n\geq 0}Q_{k,a}(n)q^n=\sum_{N_1\geq N_2\geq\cdots\geq
		N_{k-1}\geq0}\frac{q^{N_1^2+N_2^2+\cdots+N_{k-1}^2+n_2+n_4+\cdots+n_{a-3}+N_{a-1}+N_{a}
			+N_{a+1}+\cdots+N_{k-1}}
	}{(q^2;q^2)_{N_1-N_2}\cdots(q^2;q^2)_{N_{k-2}-N_{k-1}}(q^2;q^2)_{N_{k-1}}}.
\end{equation}	

\item If $a$ is even, then \begin{equation}
	\sum_{n\geq 0}Q_{k,a}(n)q^n=\sum_{N_1\geq N_2\geq\cdots\geq
		N_{k-1}\geq0}\frac{q^{N_1^2+N_2^2+\cdots+N_{k-1}^2+n_2+n_4+\cdots+n_{a-2}+N_{a}
			+N_{a+1}+\cdots+N_{k-1}}
	}{(q^2;q^2)_{N_1-N_2}\cdots(q^2;q^2)_{N_{k-2}-N_{k-1}}(q^2;q^2)_{N_{k-1}}}.
\end{equation}	
\end{itemize}
	 
\end{thm}

 In \cite{And10}, one of the open problems Andrews listed was to investigate the function
\begin{equation}\label{eqARRG}
	\sum_{N_1\geq N_2\geq \cdots \geq N_{k-1}\geq 0}\frac{q^{N_1^2+N_2^2+\cdots+N_{k-1}^2}x^{N_1+N_2+\cdots+N_{k-1}}(-yq)_{n_
			1}(-yq)_{n_2}\ldots(-yq)_{n_{k-1}}}
	{(q^2;q^2)_{N_1-N_2}\ldots(q^2;q^2)_{N_{k-2}-N_{k-1}}(q^2;q^2)_{N_{k-1}}}.
\end{equation} 

Kur\c{s}ung\"{o}z \cite{kur09}  used the Gordon marking to  interprete the left hand side  of identity \eqref{abRRG} as the generating function of $B_{k,a}(n)$. Also by this combinatrial tool he gave the answer of Andrews' question. 
To give the combinatorial interpretation of \eqref{eqARRG}, Kur\c{s}ung\"{o}z  \cite{kur10} introduced the cluster in Gordon marking.
We will describe the Gordon marking of partitions and other relative definitions which are introduced by
Kur\c{s}ung\"{o}z \cite{kur09} in Section 2. And then state the combintorial interpretation of \eqref{eqARRG} of Kur\c{s}ung\"{o}z in terms of  Gordon marking of partitions, as well as our interpretation by using lattice path.

This paper is  organized as follows. In Section 2, we shall give the definitions of Gordon marking and relative results, by which we can state the theorem of Kur\c{s}ung\"{o}z in which the combinatorial interpretation of \eqref{eqARRG} is given. Inspired by this, we introduce some definitions in lattice paths then give a combinatorial  interpretation of \eqref{eqARRG} in terms of  lattice path. In Section 3, we give the most important part of this paper, in which we  construct a bijection between partitions and lattice paths. In Section 4, in view of  the bijection in section 3, we give the proofs of  Theorem \ref{thmwp1},\ref{thmwp2},\ref{thmle} and \ref{index}.

\section{Combinatorial interpretation of Andrews' function}

In this section, we shall describe the Gordon marking of partitions and other relative definitions which are introduced by
Kur\c{s}ung\"{o}z \cite{kur09}.  Then we give the combinatorial interpretation of \eqref{eqARRG} given by  Kur\c{s}ung\"{o}z \cite{kur09} via the Gordon marking of partition. By defining some indices on lattice path we give our combinatorial  interpretation of \eqref{eqARRG} in terms of lattice path.

Firstly, let us recall the definition of Gordon marking.
\begin{defi}
	A Gordon marking of an ordinary partition $\lambda$ is an assignment of
	positive integers (marks) to parts of $\lambda$ such that any two equal parts,
	as well as
	any two nearly equal parts $j$ and $j+1$ are assigned different marks,
	and the marks are as small as possible assuming that the
	marks are assigned to the parts in increasing order.
\end{defi} For example, the Gordon marking of
\[ \lambda=(9,8,6,6,5,5,4,4,3,2,1,1)\]
can be expressed as follows
\begin{equation} \label{gm}
	\lambda=\setcounter{MaxMatrixCols}{9}\begin{bmatrix}
		\ &\ &\ &\ &5&\ &\ &\ &\ \\[3pt]
		\ &2&\ &4&\ &6 &\ &\ &\ \\[3pt]
		1&\ &\ &4&\ &6&\ &\ &9\\[3pt]
		1&\ &3&\ &5&\ &\ &8&\
	\end{bmatrix}\begin{matrix}
		4\\[3pt]3\\[3pt]2\\[3pt]1
	\end{matrix} \; ,
\end{equation}
where the marks are listed outside the brackets. The Gordon marking of a partition is  a way of representing a partition. 

To give the combinatorial interpretation of \eqref{eqARRG}, Kur\c{s}ung\"{o}z \cite{kur10} introduced the cluster in Gordon marking.
\begin{defi}
	An $r$-cluster in $\lambda=\lambda_1+\cdots+\lambda_m$ is a sub-partition $\lambda_{i_1}\leq \cdots\leq \lambda_{i_r}$ such that $\lambda_{i_j}$ is $j$-marked for $j=1,\ldots,r$, $\lambda_{i_{j+1}}-\lambda_{i_j}\leq 1$  for  $j = 1,\ldots,r-1$, and
	there are no $(r+1)$-marked parts that are equal to $\lambda_{i_r}$ or $\lambda_{i_{r}}+1$.
\end{defi}
The following proposition is obviously.
\begin{pro}
	Any partition $\lambda=\lambda_1 +\cdots+\lambda_m$ along with its Gordon marking has
	a unique decomposition into non-overlapping $r$-clusters.
\end{pro}
Kur\c{s}ung\"{o}z \cite{kur10} also proved the following proposition and gave the definitions of the parity of a cluster, lower even $r$-cluster parity index as well as the full
lower even cluster parity index of a partition. 

\begin{defi}
	Parity of an $r$-cluster is the opposite parity of the number of even parts in that $r$-cluster. Lower even $r$-cluster parity index of a partition $\lambda$ is the number of
	times that the $r$-cluster parity changes from the $r$-cluster with the smallest $r$-marked
	part to the one with the largest $r$-marked part, beginning with an even $r$-cluster parity.
\end{defi}
\begin{defi}
	Given a partition $\lambda=\lambda_1+\cdots+\lambda_m$ enumerated by $B_{k,k}(m,n)$, the full	lower even cluster parity index of $\lambda$ is the sum of all lower even $1-, 2-,\ldots,(k-1)-$
	cluster parity indices.
\end{defi}

We rewrite  Kur\c{s}ung\"{o}z's result as follows which gives a combintorial interpretation of \eqref{eqARRG} in terms of  partition.

\begin{thm}[Kur\c{s}ung\"{o}z \cite{kur10}]\label{kur}	Let $B_{k,k}(l,m,n)$ denote the number of partitions  enumerated by $B_{k,k}(m,n)$ whose  full lower even cluster parity index are $l$, then the generating function of $B_{k,k}(l,m,n)$ is \eqref{eqARRG}, that is,
	\begin{align*}
		&\sum_{l,m,n\geq0}B_{k,k}(l,m,n)y^lx^mq^n\\&\qquad=	\sum_{N_1\geq N_2\geq \cdots \geq N_{k-1}\geq 0}\frac{q^{N_1^2+N_2^2+\cdots+N_{k-1}^2}x^{N_1+N_2+\cdots+N_{k-1}}(-yq)_{n_1}(-yq)_{n_2}\ldots(-yq)_{n_{k-1}}}
		{(q^2;q^2)_{N_1-N_2}\ldots(q^2;q^2)_{N_{k-2}-N_{k-1}}(q^2;q^2)_{N_{k-1}}}.
	\end{align*} 
	
\end{thm}

In this paper, we also give an combinatorial interpretation of \eqref{eqARRG} in terms of  lattice path. We need some definitions.
\begin{defi}\label{defiparitypath}
	If the weight of the peak is $x$ and the relative height of the peak is $r$. Then, the parity of the peak  defined as the opposite parity of $x-r$.  Lower even $r$-peak parity index of a lattice path $L$ is the number of
	times that the parity of peaks with  relative height $r$  changes from the leftmost to the rightmost, beginning with a peak with even parity and relative hight $r$.
\end{defi}

\begin{defi}
	Given a lattice path $L$ satisfies the special $(k,k)$-conditions, the full
	lower even peak parity index of $L$ is the sum of all lower even $1-, 2-,\ldots,(k-1)-$ 
	peaks parity indices.
\end{defi}
Then we are ready to give our result.
\begin{thm}\label{index}
	Let $E_{k,k}(l,m,n)$ denote the number of lattice paths  enumerated by $E_{k,k}(n)$ in which the  sum of  relative height of all peaks is $m$ and full lower even peak parity index is $l$, then
	\[E_{k,k}(l,m,n)=B_{k,k}(l,m,n).\]
\end{thm}
That is to say, the generating function of $E_{k,k}(l,m,n)$ is also \eqref{eqARRG}.

\section{The bijection between Gordon marking of a partition and a lattice path}

In this section, we shall give a bijection between partition $\lambda_0$ enumerated by $B_{k,a}(n)$ and a lattice path $L$ enumerated by $E_{k,a}(n)$.

To give the bijection between Gordon marking of partitions and lattice paths, firstly, we shall describe how to  constructe a partition in terms of clusters in the Gordon marking. 
 Let $\mathcal{B}_{k,a}(n)$ denote the set of partitions enumerated by $B_{k,a}(n)$, and $\mathcal{B}_{N_1,\ldots,N_{k-1};a}(n)$ be the subset of $\mathcal{B}_{k,a}(n)$ in which the partitions have $N_i$ $i$-marked parts. Let 
 \[\mathcal{B}_{k,a}=\bigcup_{n\geq0}\mathcal{B}_{k,a}(n),\]
 \[\mathcal{B}_{N_1,\ldots,N_{k-1};a}=\bigcup_{n\geq0}\mathcal{B}_{N_1,\ldots,N_{k-1};a}(n).\] 
  Given a partition $\lambda_0\in\mathcal{B}_{N_1,\ldots,N_{k-1};a}(n)$,  we firstly relate $\lambda$ with a partition $\mu\in \mathcal{B}_{N_1,\ldots,N_{k-1};a}$ and $k-1$ partitions.
 
In the Gordon marking of $\lambda_0$, there are  $N_r$ $r$-marked parts and $n_r$ $r$-cluster, where $N_r=n_r+\cdots+n_{k-1}$, for  $r=1,2,\ldots, k-1$. We will map $\lambda_0$ to a partition $\mu\in\mathcal{B}_{N_1,\ldots,N_{k-1};a}$ who has smallest weight with $n_r$ $r$-cluster, for $r=1,2,\ldots,k-1$, and $k-1$ partitions $\pi^{(1)},\pi^{(2)},\cdots,\pi^{(k-1)}$ where $\pi^{(i)}$ has $n_i$ nonnegative parts for $i=1,\ldots,  k-1$. 

We sketch the map $\phi$ by using the backward move defined by   Kur\c{s}ung\"{o}z \cite{kur09}.

\noindent $\phi$: From $\lambda_0$ to $\mu$ and $\pi^{(1)},\pi^{(2)},\cdots,\pi^{(k-1)}$.

 We start at  $\lambda=\lambda_0$  and $\mu=\emptyset$, $\pi^{(1)}=\pi^{(2)}=\cdots=\pi^{(k-1)}=\emptyset$.
For $r$ from $k-1$ to $1$.
\begin{itemize}
\item [Step 1.]	There are $m$ clusters in $\mu$ and $r$ is the largest mark in $\lambda$. We choose the $r$-cluster with the smallest $r$-marked part, and denote it by $\lambda_0^{(r)}$. We employ the backward move of the $r$-th kind to $\lambda^{(r)}_0$ until the cluster that contains $\lambda^{(r)}_0$ have $a-1$ parts equal to $2m+1$ and $r-a+1$ parts equal to $2m+2$, if $r\geq a$; or contains $r$ parts equal to $2m+1$, if $r\leq a-1$.

\item[Step 2.] Move the cluster which contains $\lambda^{(r)}_0$ to $\mu$. If there are no $r$-marked part in $\lambda$ and $r\geq 2$, let $r=r-1$ and back to Step 1; if there are some $r$-marked parts in $\lambda$ back to Step 1 directly; if $\lambda=\emptyset$, go to Step 3.

\item[Step 3.] We get a Gordon marking of partition $\mu$ with $N_1+\cdots+N_{k-1}$ parts, $N_1$ cluster and $k-1$ partitions $\pi^{(1)},\pi^{(2)},\cdots,\pi^{(k-1)}$ where $\pi^{(r)}$ has $n_r$  non-negative parts, $r= 1,\ldots,k-1$.  We can also derive that $|\lambda_0|=|\mu|+|\pi^{(1)}|+|\pi^{(2)}|+\cdots+|\pi^{(k-1)}|$.\qed
\end{itemize}

We shall give an example to show this map.
\begin{exam}\label{exam2}Let $k=4,\ a=3$, and  partition $\lambda_0=(13,11,11,11,9,8,6,6,5,4,3,3,2,1)\in\mathcal{B}_{4,3} $.  The Gordon marking of $\lambda_0$ is displaied as follows.
\begin{equation*}\label{lambda}
	\lambda_0=\setcounter{MaxMatrixCols}{13}\begin{bmatrix}
		\ &\ &3&\ &\ &6&\ &\ &\ &\ &11&\ &\ \\[3pt]
		\ &2&\ &4&\ &6&\ &\ &9&\ &11&\ &\ \\[3pt]
		1&\ &3&\ &5&\ &\ &8&\ &\ &11&\ &13
	\end{bmatrix}\begin{matrix}
		3\\[3pt]2\\[3pt]1
	\end{matrix}\;,
\end{equation*}
Accordingly, $\lambda_0$ has three $3$-cluster, two $2$-cluster and one $1$-cluster and  $|\lambda_0|=93$.
We begain with the $3$-cluster $(1_1,2_2,3_3)$ which has  the smallest $3$-marked part $3$.
We employ the backward move to the $3$-marked $3$ two times as follows:
\begin{equation*}\label{lambda}
	\lambda_0\rightarrow \setcounter{MaxMatrixCols}{13}\begin{bmatrix}
		\ &\mathbf{2}&\ &\ &\ &6&\ &\ &\ &\ &11&\ &\ \\[3pt]
		\ &\mathbf{2}&\ &4&\ &6&\ &\ &9&\ &11&\ &\ \\[3pt]
	\mathbf{1}&\ &3&\ &5&\ &\ &8&\ &\ &11&\ &13
	\end{bmatrix}\begin{matrix}
		3\\[3pt]2\\[3pt]1
	\end{matrix}\;\rightarrow \setcounter{MaxMatrixCols}{13}\begin{bmatrix}
	\ &\mathbf{2}&\ &\ &\ &6&\ &\ &\ &\ &11&\ &\ \\[3pt]
	\mathbf{1}&\ &\ &4&\ &6&\ &\ &9&\ &11&\ &\ \\[3pt]
	\mathbf{1}&\ &3&\ &5&\ &\ &8&\ &\ &11&\ &13
\end{bmatrix}\begin{matrix}
3\\[3pt]2\\[3pt]1
\end{matrix}\;,
\end{equation*}
We move the  $3$-cluster having $3$-marked $2$ to $\mu$ to get $\mu=(1_1,1_2,2_3)$, $\pi^{(3)}=(2)$ and \begin{equation*} \lambda=
\setcounter{MaxMatrixCols}{13}\begin{bmatrix}
	\ &\ &\ &\ &\ &\mathbf{6}&\ &\ &\ &\ &11&\ &\ \\[3pt]
\ &\ &\ &4&\ &\mathbf{6}&\ &\ &9&\ &11&\ &\ \\[3pt]
	\ &\ &3&\ &\mathbf{5}&\ &\ &8&\ &\ &11&\ &13
\end{bmatrix}\begin{matrix}
	3\\[3pt]2\\[3pt]1
\end{matrix}\;,
\end{equation*}
In view of  the smallest $3$-marked part is $6$, we consider the $3$-cluster $(5_1,\ 6_2,\ 6_3)$. We repeat $3$ times backward move to $3$-marked $6$ to make the $3$-cluster $(5,6,6)$ to be $3$-cluster $(3,3,4)$ and move this cluster to $\mu$ as well as add $3$ to $\pi^{(3)}$, the procedures are showed below, 
\begin{equation*} 
\lambda=	\setcounter{MaxMatrixCols}{13}\begin{bmatrix}
	\ &\ &\ &\mathbf{4}&\ &\ &\ &\ &\ &\ &11&\ &\ \\[3pt]
	\ &\ &\mathbf{3}&\ &\ &6&\ &\ &9&\ &11&\ &\ \\[3pt]
	\ &\ &\mathbf{3}&\ &5&\ &\ &8&\ &\ &11&\ &13
\end{bmatrix}\begin{matrix}
	3\\[3pt]2\\[3pt]1
\end{matrix}\;, 
\end{equation*}
and $\mu$ is as follows.
\begin{equation*}\lambda=	\setcounter{MaxMatrixCols}{13}\begin{bmatrix}
		\ &\ &\ &\ &\ &\ &\ &\ &\ &\ &\mathbf{11}&\ &\ \\[3pt]
		\ &\ &\ &\ &\ &6&\ &\ &9&\ &\mathbf{11}&\ &\ \\[3pt]
		\ &\ &\ &\ &5&\ &\ &8&\ &\ &\mathbf{11}&\ &13
	\end{bmatrix}\begin{matrix}
		3\\[3pt]2\\[3pt]1
	\end{matrix}\;,\qquad \mu=\setcounter{MaxMatrixCols}{5}\begin{bmatrix}
	\ &2&\ &4&\ \\[3pt]
	1&\ &3&\ &\ \\[3pt]
	1&\ &3&\ &\ 
\end{bmatrix}\begin{matrix}
	3\\[3pt]2\\[3pt]1
\end{matrix}\;.
\end{equation*}
Now the last $3$-cluster is  $(11_1,11_2,11_3)$, we also employ the backward move $8$ times to the $3$-marked $11$, then we get a $3$-marked cluster $(5,6,6)$ and $\pi^{(3)}=(8,3,2)$.
\begin{equation*} 
	\lambda=	\setcounter{MaxMatrixCols}{13}\begin{bmatrix}
		\ &\ &\ &\ &\ &\mathbf{6}&\ &\ &\ &\ &\ &\ &\ \\[3pt]
		\ &\ &\ &\ &\ &\mathbf{6}&\ &8 &\ &\ &11&\ &\ \\[3pt]
		\ &\ &\ &\ &\mathbf{5} &\ &7&\ &\ &\ &11&\ &13
	\end{bmatrix}\;.
\end{equation*}

 Moving this cluster to $\mu$ we get $\lambda$ and $\mu$ as follows: 
\begin{equation*} 
	\lambda=	\setcounter{MaxMatrixCols}{13}\begin{bmatrix}
		\ &\ &\ &\ &\ &\ &\ &\ &\ &\ &\ &\ &\ \\[3pt]
		\ &\ &\ &\ &\ &\ &\ &\mathbf{8} &\ &\ &11&\ &\ \\[3pt]
		\ &\ &\ &\ &\ &\ &\mathbf{7}&\ &\ &\ &11&\ &13
	\end{bmatrix}\begin{matrix}
		3\\[3pt]2\\[3pt]1
	\end{matrix}\;,\qquad \mu=\setcounter{MaxMatrixCols}{6}\begin{bmatrix}
		\ &2&\ &4&\ &6\\[3pt]
		1&\ &3&\ &5&\ \\[3pt]
		1&\ &3&\ &5&\ 
	\end{bmatrix}\begin{matrix}
		3\\[3pt]2\\[3pt]1
	\end{matrix}\;.
\end{equation*}

Now there is no $3$-cluster in $\lambda$, let $r=2$ back to Step 1. We make backward move to the smallest $2$-marked part $8$ for one time, then get a $2$-cluster $(7_1,7_2)$ and $\pi^{(2)}=(1)$, after move this cluster to $\mu$, we get $\lambda$ and $\mu$ as follows:
\begin{equation*} 
	\lambda=	\setcounter{MaxMatrixCols}{12}\begin{bmatrix}
		\ &\ &\ &\ &\ &\ &\ &\ &\ &\ &\ &\ \\[3pt]
		\ &\ &\ &\ &\ &\ &\ &\ &\ &\mathbf{11}&\ &\ \\[3pt]
		\ &\ &\ &\ &\ &\ &\ &\ &\ &\mathbf{11}&\ &13
	\end{bmatrix}\begin{matrix}
		3\\[3pt]2\\[3pt]1
	\end{matrix}\;,\qquad \mu=\setcounter{MaxMatrixCols}{8}\begin{bmatrix}
		\ &2&\ &4&\ &6&\ &\ \\[3pt]
		1&\ &3&\ &5&\ &7&\ \\[3pt]
		1&\ &3&\ &5&\ &7
	\end{bmatrix}\begin{matrix}
		3\\[3pt]2\\[3pt]1
	\end{matrix}\;.
\end{equation*}
In view of the $2$-cluster $(11_1,11_2)$ we make the backward move to $2$-marked $11$ four times to get a $2$-cluster $(9_1,9_2)$ and $\pi^{(2)}=(4,2)$. After moving this cluster to $\mu$, we get $\lambda$ and $\mu$ as follows:
\begin{equation*} 
	\lambda=	\setcounter{MaxMatrixCols}{10}\begin{bmatrix}
		\ &\ &\ &\ &\ &\ &\ &\ &\ \\[3pt]
		\ &\ &\ &\ &\ &\ &\ &\ &\ \\[3pt]
		\ &\ &\ &\ &\ &\ &\ &\ &\mathbf{13}
	\end{bmatrix}\begin{matrix}
		3\\[3pt]2\\[3pt]1
	\end{matrix}\;,\qquad \mu=\setcounter{MaxMatrixCols}{9}\begin{bmatrix}
		\ &2&\ &4&\ &6&\ &\ &\ \\[3pt]
		1&\ &3&\ &5&\ &7&\ &9\\[3pt]
		1&\ &3&\ &5&\ &7&\ &9
	\end{bmatrix}\begin{matrix}
		3\\[3pt]2\\[3pt]1
	\end{matrix}\;.
\end{equation*}
There is only one $1$-cluster $(13_1)$ we employ the backward move twice to get a cluster $(11_1)$ and $\pi^{(1)}=(2)$. Moving this cluster to $\mu$ we get $\lambda=\emptyset$ and $\mu$ as follows:
\begin{equation*} \mu=\setcounter{MaxMatrixCols}{11}\begin{bmatrix}
	\ &2&\ &4&\ &6&\ &\ &\ &\ &\  \\[3pt]
	1&\ &3&\ &5&\ &7&\ &9&\ &\ \\[3pt]
	1&\ &3&\ &5&\ &7&\ &9&\ &11
\end{bmatrix}\begin{matrix}
	3\\[3pt]2\\[3pt]1
\end{matrix}\;.
\end{equation*}
The three partitions are $\pi^{(1)}=(2),\pi^{(2)}=(4,1)$ and $\pi^{(3)}=(8,3,2)$. We can check that $|\mu|=73$, $|\lambda_0|=|\mu|+|\pi^{(1)}|+|\pi^{(2)}|+|\pi^{(3)}|$.
\end{exam}

Now we proceed to  construct a lattice path $L$ from partition $\mu$ and $\pi^{(1)},\ldots,\pi^{(k-1)}$.  The procedures of the construction  will show  us that there is a bijection between the partitions in  $\mathcal{B}_{N_1,\ldots,N_{k-1};a}(n)$ and the lattice paths in $\mathcal{E}_{N_1,\ldots,N_{k-1};a}(n)$. We also show that $L$ has $N_1$ peaks and in which $n_r$ peaks have relative height $r$. 

We shall use the ``volcanic uplift"  which introduced by Bressoud in \cite{bre87}.

\begin{defi} The volcanic uplift  performs as following. At each peak we break the lattice path,
spread it apart by two units, and insert a new peak whose height is one more
than the old. 
\end{defi}We give the following proposition.

\begin{pro} The volcanic uplift changes the parity of the weight of each peak in the lattice path which it perfomed on.
\end{pro}
\pf The peaks increase in weight by $1, 3, 5, \ldots$ successively. \qed

Now we give the map $\psi$:

\noindent
$\psi$: For $i$ from $k-1$ to $1$.
\begin{itemize}
	\item [Step 1] If there are $n_i$ $i$-cluster in $\mu$,   we  insert $n_i$ peaks with relative height  one  at the start point of $L$. More precisely, if $i\geq a$,  the start point is $(0,k-i-1)$, hence the peaks at $(1,k-i),\ (3,k-i),\ldots (2n_i-1,k-i)$; if  $i\leq a-1$,  the start  point is $(0,k-a)$, then the peaks at $(1,k-a+1),\ (3,k-a+1),\ \ldots,(2n_i-1,k-a+1)$.
	
\item[Step 2] If $i\geq a$, we add an SE step start at $(0,k-i)$, hence the weight of  each peak enlarge one.

\item[Step 3]	For $\pi^{(i)}$ is a partition with $n_i$ non-negative parts, then for $j=1,\ldots, n_{i}$, from right to left we let the $j$-th peak of relative height one increase it weight $\pi^{(i)}_j$ by move the peak to the right $\pi^{(i)}_j$.  The way of moving the peak of relative height one to right is the same of that introduced by Bressoud \cite{bre87}.
	
\item[Step 4] If $i\geq 2$, we want to increase the height of each peak, in doing so, we perform  the  ``volcanic uplift"  and then let $i=i-1$, back to Step 1. If $i=1$, we get the lattice path $L$.
\end{itemize}

The procedures of construct lattice path $L$ show us that $L$  satisfies the special $(k,a)$-conditions having $n_r$ peaks with relative height $r$. One can easily check that the major index of $L$  equals $|\mu|+|\pi^{(1)}|+|\pi^{(2)}|+\cdots+|\pi^{(k-1)}|=|\lambda_0|$.

Let $\mathcal{E}_{k,a}(n)$ denote the set of lattice paths enumerated by $E_{k,a}(n)$, and $\mathcal{E}_{N_1,\ldots,N_{k-1};a}(n)$ be the subset of $\mathcal{E}_{k,a}(n)$ in which the lattice paths  have $N_i$ peaks with relative height greater than or equal to  $i$. Let 
\[\mathcal{E}_{k,a}=\bigcup_{n\geq0}\mathcal{E}_{k,a}(n),\]
\[\mathcal{E}_{N_1,\ldots,N_{k-1};a}=\bigcup_{n\geq0}\mathcal{E}_{N_1,\ldots,N_{k-1};a}(n).\] 

We have the observations that $\phi\times \psi$ is  a bijection between $\mathcal{B}_{N_1,\ldots,N_{k-1};a}(n)$ and $\mathcal{E}_{N_1,\ldots,N_{k-1};a}(n)$.  

\begin{exam}We shall construct a lattice path satisfies the special $(4,3)$-conditions from $\mu$ and partitions $\pi^{(1)},\pi^{(2)},\pi^{(3)}$ given in Example \ref{exam2}.
Firstly, we consider the three  $3$-clusters in $\mu$, which generate three peaks with relative height one at $(1,1),(3,1),(5,1)$.
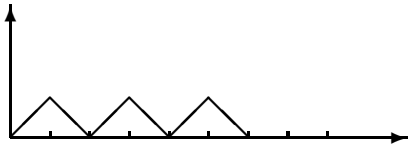
\begin{center}
\begin{picture}(150,70)
\thicklines
\multiput(25,10)(15,0){8}{\line(0,1){2}}
\put(10,10){\vector(1,0){150}}
\put(10,10){\vector(0,1){50}}
\put(10,10){\line(1,1){15}}
\put(25,25){\line(1,-1){15}}
\put(40,10){\line(1,1){15}}
\put(55,25){\line(1,-1){15}}
\put(70,10){\line(1,1){15}}
\put(85,25){\line(1,-1){15}}
\end{picture}
\end{center}
Due to the fact that  $3=i\geq a=3$, we should add an SE step at $(0,1)$,
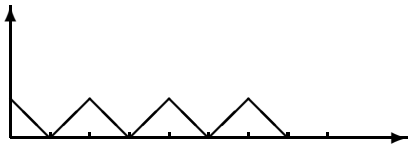
\begin{center}
	\begin{picture}(150,70)
		\thicklines
		\multiput(25,10)(15,0){8}{\line(0,1){2}}
		\put(10,10){\vector(1,0){150}}
		\put(10,10){\vector(0,1){50}}
		\put(10,25){\line(1,-1){15}}
		\put(25,10){\line(1,1){15}}
		\put(40,25){\line(1,-1){15}}
		\put(55,10){\line(1,1){15}}
		\put(70,25){\line(1,-1){15}}
		\put(85,10){\line(1,1){15}}
		\put(100,25){\line(1,-1){15}}
	\end{picture}
\end{center}
In view of $\pi_3=(8,3,2)$, we  move  the  rightmost peak to right for $8$ times,  the second peak for $3$ times and the leftmost peak for twice according to the  rules Bressound introduced in \cite{Bre80}, hence we get the following lattice path.
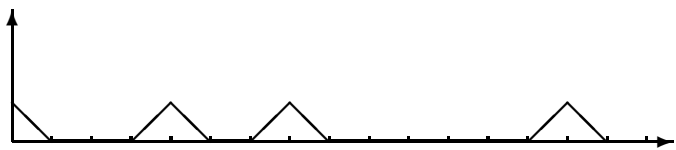
\begin{center}
	\begin{picture}(300,70)
		\thicklines
		\multiput(25,10)(15,0){16}{\line(0,1){2}}
		\put(10,10){\vector(1,0){250}}
		\put(10,10){\vector(0,1){50}}
		\put(10,25){\line(1,-1){15}}
		\put(25,10.5){\line(1,0){30}}
		\put(55,10){\line(1,1){15}}
		\put(70,25){\line(1,-1){15}}
		\put(85,10.5){\line(1,0){15}}
		\put(100,10){\line(1,1){15}}
		\put(115,25){\line(1,-1){15}}
		\put(130,10.5){\line(1,0){75}}
		\put(205,10){\line(1,1){15}}
		\put(220,25){\line(1,-1){15}}
	\end{picture}
\end{center}
Now we  perform  the volcanic uplift and let $i=2$.
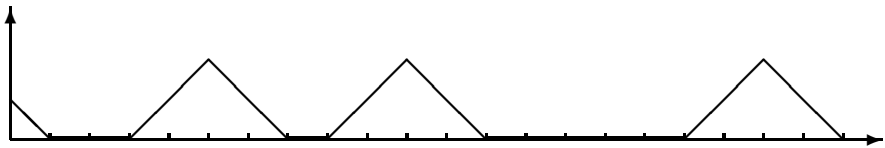
\begin{center}
\begin{picture}(340,70)
	\thicklines
	\multiput(25,10)(15,0){21}{\line(0,1){2}}
	\put(10,10){\vector(1,0){330}}
	\put(10,10){\vector(0,1){50}}
	\put(10,25){\line(1,-1){15}}
	\put(25,10.5){\line(1,0){30}}
	\put(55,10){\line(1,1){30}}
	\put(85,40){\line(1,-1){30}}
	\put(115,10.5){\line(1,0){15}}
	\put(130,10){\line(1,1){30}}
	\put(160,40){\line(1,-1){30}}
	\put(190,10.5){\line(1,0){75}}
	\put(265,10){\line(1,1){30}}
	\put(295,40){\line(1,-1){30}}
\end{picture}
\end{center}
There are two $2$-cluster in $\mu$. We insert two peaks with relative height one.
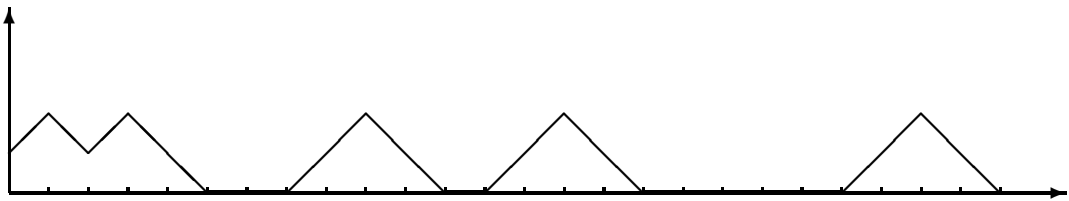
\begin{center}
	\begin{picture}(400,80)
		\thicklines
		\multiput(25,10)(15,0){25}{\line(0,1){2}}
		\put(10,10){\vector(1,0){400}}
		\put(10,10){\vector(0,1){70}}
		\put(10,25){\line(1,1){15}}
		\put(25,40){\line(1,-1){15}}
		\put(40,25){\line(1,1){15}}
		\put(55,40){\line(1,-1){15}}
		\put(70,25){\line(1,-1){15}}
		\put(85,10.5){\line(1,0){30}}
		\put(115,10){\line(1,1){30}}
		\put(145,40){\line(1,-1){30}}
		\put(175,10.5){\line(1,0){15}}
		\put(190,10){\line(1,1){30}}
		\put(220,40){\line(1,-1){30}}
		\put(250,10.5){\line(1,0){75}}
		\put(325,10){\line(1,1){30}}
		\put(355,40){\line(1,-1){30}}
	\end{picture}
\end{center}
Because of $\pi_2=(4,1)$,  we firstly move the rightmost  peak who has  relative height one  to right $4$ times to get the following lattice path:
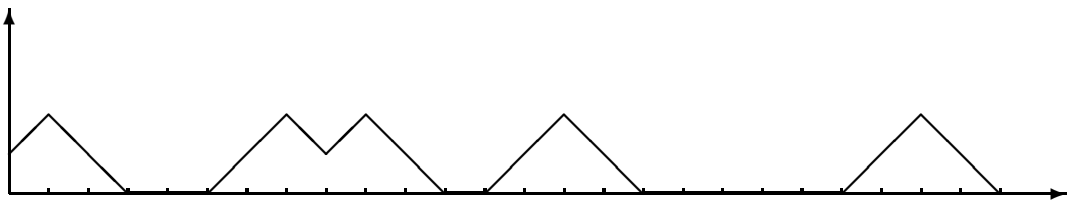
\begin{center}
	\begin{picture}(400,80)
		\thicklines
		\multiput(25,10)(15,0){25}{\line(0,1){2}}
		\put(10,10){\vector(1,0){400}}
		\put(10,10){\vector(0,1){70}}
		\put(10,25){\line(1,1){15}}
		\put(25,40){\line(1,-1){15}}
		\put(130,25){\line(1,1){15}}
		\put(115,40){\line(1,-1){15}}
		\put(40,25){\line(1,-1){15}}
		\put(55,10.5){\line(1,0){30}}
		\put(85,10){\line(1,1){30}}
		\put(145,40){\line(1,-1){30}}
		\put(175,10.5){\line(1,0){15}}
		\put(190,10){\line(1,1){30}}
		\put(220,40){\line(1,-1){30}}
		\put(250,10.5){\line(1,0){75}}
		\put(325,10){\line(1,1){30}}
		\put(355,40){\line(1,-1){30}}
	\end{picture}
\end{center}
Then we move the leftmost  peak with  relative height one  to right once to get the following lattice path:
\begin{center}
	\begin{picture}(400,80)
		\thicklines
		\multiput(25,10)(15,0){25}{\line(0,1){2}}
		\put(10,10){\vector(1,0){400}}
		\put(10,10){\vector(0,1){70}}
		\put(25,10){\line(1,1){15}}
		\put(10,25){\line(1,-1){15}}
		\put(130,25){\line(1,1){15}}
		\put(115,40){\line(1,-1){15}}
		\put(40,25){\line(1,-1){15}}
		\put(55,10.5){\line(1,0){30}}
		\put(85,10){\line(1,1){30}}
		\put(145,40){\line(1,-1){30}}
		\put(175,10.5){\line(1,0){15}}
		\put(190,10){\line(1,1){30}}
		\put(220,40){\line(1,-1){30}}
		\put(250,10.5){\line(1,0){75}}
		\put(325,10){\line(1,1){30}}
		\put(355,40){\line(1,-1){30}}
	\end{picture}
\end{center}

Now we perform the  volcanic uplift and let $i=1$.
\begin{center}
	\begin{picture}(400,100)
		\thicklines
		\multiput(20,10)(10,0){35}{\line(0,1){2}}
		\put(10,10){\vector(1,0){380}}
		\put(10,10){\vector(0,1){70}}
		\put(20,10){\line(1,1){20}}
		\put(10,20){\line(1,-1){10}}
		
		\put(40,30){\line(1,-1){20}}
		\put(60,10.5){\line(1,0){20}}
		\put(80,10){\line(1,1){30}}
		\put(110,40){\line(1,-1){20}}\put(130,20){\line(1,1){20}}
		\put(150,40){\line(1,-1){30}}
		\put(180,10.5){\line(1,0){10}}
		\put(190,10){\line(1,1){30}}
		\put(220,40){\line(1,-1){30}}
		\put(250,10.5){\line(1,0){50}}
		\put(300,10){\line(1,1){30}}
		\put(330,40){\line(1,-1){30}}
	\end{picture}
\end{center}
There is one $1$-cluster in $\mu$, so we add one peak with relative height one to the lattice.
\begin{center}
	\begin{picture}(400,100)
		\thicklines
		\multiput(20,10)(10,0){40}{\line(0,1){2}}
		\put(10,10){\vector(1,0){420}}
		\put(10,10){\vector(0,1){70}}
		
		\put(10,20){\line(1,1){10}}
		\put(20,30){\line(1,-1){10}}
		
		\put(30,20){\line(1,-1){10}}
		\put(40,10){\line(1,1){20}}
		\put(60,30){\line(1,-1){20}}
		\put(80,10.5){\line(1,0){20}}
		\put(100,10){\line(1,1){30}}
		\put(130,40){\line(1,-1){20}}\put(150,20){\line(1,1){20}}
		\put(170,40){\line(1,-1){30}}
		\put(200,10.5){\line(1,0){10}}
		\put(210,10){\line(1,1){30}}
		\put(240,40){\line(1,-1){30}}
		\put(270,10.5){\line(1,0){50}}
		\put(320,10){\line(1,1){30}}
		\put(350,40){\line(1,-1){30}}
	\end{picture}
\end{center}
Considering that $\pi^{(1)}=(2)$, we move the peak with relative one to the right twice.
\begin{center}
	\begin{picture}(400,100)
		\thicklines
		\multiput(20,10)(10,0){40}{\line(0,1){2}}
		\put(10,10){\vector(1,0){420}}
		\put(10,10){\vector(0,1){70}}
		
		\put(20,10){\line(1,1){20}}
		\put(10,20){\line(1,-1){10}}
		
		\put(40,30){\line(1,-1){10}}
		\put(50,20){\line(1,1){10}}
		\put(60,30){\line(1,-1){20}}
		\put(80,10.5){\line(1,0){20}}
		\put(100,10){\line(1,1){30}}
		\put(130,40){\line(1,-1){20}}\put(150,20){\line(1,1){20}}
		\put(170,40){\line(1,-1){30}}
		\put(200,10.5){\line(1,0){10}}
		\put(210,10){\line(1,1){30}}
		\put(240,40){\line(1,-1){30}}
		\put(270,10.5){\line(1,0){50}}
		\put(320,10){\line(1,1){30}}
		\put(350,40){\line(1,-1){30}}
	\end{picture}
\end{center}
At last, we get a lattice path with one peak whose relative height is one and weight is $5$, two peaks whose relative height are two and weight are $3$ and $16$ as well as three peaks whose relative height are $3$ and weights are $12$, $23$ and $34$. The major index of this lattice path is $5+3+16+12+23+34=93=|\lambda_0|$. 
\end{exam}

\section{The proofs  of theorems}

In this section, we shall give the proofs of Theorem \ref{thmwp1}, \ref{thmwp2}, \ref{thmle} and \ref{index}. 
Let $\mu$ denote the partition in $\mathcal{B}_{N_1,\ldots,N_{k-1};a}$  with smallest weight. Then we know that the $r$-marked parts in $\mu$ are $1,3,\ldots, 2N_r-1$, if $r\leq a-1$, or $2,4,\ldots, 2N_r$, if $r\geq a$.

\noindent{\bf The proof of Theorem \ref{thmwp1}.}
Let $\mathcal{W}_{N_1,\ldots,N_{k-1};a}$ be the subset of $\mathcal{B}_{N_1,\ldots,N_{k-1};a}$ in which the partitions satisfy that even parts appear even times. For a partition $\lambda\in \mathcal{W}_{N_1,\ldots,N_{k-1};a}$, Kur\c{s}ung\"{o}z \cite{kur09} has proved that $\lambda$ can be constructed from $\mu$ by a series forward move. We  describe these procedures in terms of clusters in Gordon marking. That is, if $r\geq a$ and $r-(a-1)$  is odd, do the  $r$-th kind forward move odd number of times  to each $n_r$ $r$-cluster, or do even number of times to each $n_r$ $r$-cluster, otherwise, for $r=1,\ldots,k-1$. Then we can get a partition in $\mathcal{W}_{N_1,\ldots,N_{k-1};a}$.  By the description in Section 2, we know that $\lambda$ corresponds to  $\mu$ and $k-1$ partitions $\pi^{(1)},\ldots,\pi^{(k-1)}$. To ensure $\lambda$ is in $\mathcal{W}_{N_1,\ldots,N_{k-1};a}$, $\pi^{(r)}$ should be  a partition with $n_r$ odd parts if $r\geq a$ and $r\equiv a\pmod{2}$, or $\pi^{(r)}$ is a partition with $n_r$ even parts, otherwise.

By the bijection between  partitions and   lattice paths given in Section 2, we shall consider what conditions the lattice path  corresponds to a partition  $\lambda \in\mathcal{W}_{N_1,\ldots,N_{k-1};a}$ should satisfy.  Let $L$ be the lattice path corresponding to $\lambda$, the bijection tells us  that for $k-1\geq r\geq 1$, there are $n_r$  peaks in $L$ with relative height $r$. As the procedures we described in section 2, we construct them by add $n_r$ peaks at beginning with relative height one.  More precisely, if $r\geq a$,  the start point is $(0,k-r-1)$, then the peaks are at $(1,k-r),\ (3,k-r),\ldots (2n_r-1,k-r)$; if  $r\leq a-1$,  the start  point is $(0,k-a)$, then the peaks are at $(1,k-a+1),\ (3,k-a+1),\ \ldots,(2n_r-1,k-a+1)$. Therefore the peaks have odd weights. 

\begin{itemize}
\item If $r\geq a$ and $r\equiv a\pmod{2}$, we need to insert an SE step at $(0,k-r)$, which makes the weight of  peaks increase one, then the peaks are all with even weights. In view of $r\equiv a\pmod{2}$, the parts in $\pi^{(r)}$ are odd, so we increase the weight odd number of times, which makes the weights of the peaks are odd with relative height one. Then after $r-a$ times insertion of SE step and  $r-1$ times volcanic uplift, the $n_r$ peaks are peaks with relative height $r$ and weights have the same parity with $r$ also the  same parity with $a$.
\item If $r\geq a$ and $r\not\equiv a\pmod{2}$, we insert an SE step at $(0,k-r)$, which makes the weight of  peaks increase one, then the peaks are all with even weight. Because of  $r\not\equiv a\pmod{2}$, the parts in $\pi^{(r)}$ are even, so we increase the weight even number of times, which makes the weight of the peaks are even with relative height one. Then after $r-a$ times  insertion of SE step and $r-1$ times volcanic uplift, the $n_r$ peaks are peaks with relative height $r$ and weights have the opposite parity with $a$, but same parity with $r$.

\item If $r\leq a-1$, the parts in $\pi^{(r)}$ are even, thus we increase  weight of each peak even number of times, which makes the peaks with relative height one have odd weight. As a result of  $r-1$ times volcanic uplift, the $n_r$ peaks are peaks with relative height $r$ and weight have the same parity with $r$.
\end{itemize}
Accordingly, we have completed the proof of Theorem \ref{thmwp1}. By Definition \ref{defiparitypath}, one can also know that all peaks are have odd parity. \qed

\noindent{\bf The proof of Theorem \ref{thmwp2}.} Let $\overline{\mathcal{W}}_{N_1,\ldots,N_{k-1};a}$ be the subset of $\mathcal{B}_{N_1,\ldots,N_{k-1};a}$ in which the partitions satisfy that odd parts appear even times.  
Given any  partition $\lambda\in \overline{\mathcal{W}}_{N_1,\ldots,N_{k-1};a}$ , Kur\c{s}ung\"{o}z \cite{kur09} has also  proved that $\lambda$ can be constructed from $\mu$ by a series forward move. We also  describe these procedures in terms of cluster in Gordon marking.  For $r$ from $1$ to $k-1$, the forward move apply to the clusters from right to left follows the following rules.
\begin{itemize}
\item 
 If $r\geq a$ and $a$  is even, we know that each odd part appears an odd number of times and there are only $l$ and $l+1$ two different parts in  each $r$-cluster in $\mu$, therefore  we do  $r$-th kind forward move  once to each  $r$-cluster from right to left, to make the clusters have even number of odd parts. If  $r\geq a$ and $a$  is odd, then each $r$-cluster in $\mu$ have even number of odd parts, then we do nothing to the $r$-cluster. 
\item If $r\leq a-1$, and $r$ is odd then there are odd number of odd parts in this cluster, then we should make the $r$-th kind forward move once to each $r$-cluster.  
\item If $r\leq a-1$, and $r$ is even then there are even number of odd parts in this cluster, then we do nothing. 

By now we have a partition $\tilde{\mu}\in\overline{\mathcal{W}}_{N_1,\ldots,N_{k-1};a}$, which is the partition in $\overline{\mathcal{W}}_{N_1,\ldots,N_{k-1};a}$ with the smallest weight. 
\end{itemize}

Then any  partition $\lambda \in\overline{\mathcal{W}}_{N_1,\ldots,N_{k-1};a}$ can be get by applying even number of times $r$-th kind forward move to each $r$-cluster from right to left successively in $\tilde{\mu}$, for $r=1,\ldots,k-1$.

By the bijection we constructed in Section 2, we can derive that $\lambda$ corresponds to a partition $\mu\in \mathcal{B}_{N_1,\ldots,N_{k-1};a}$ and $k-1$ partitions $\pi^{(1)},\pi^{(2)},\ldots,\pi^{(k-1)}$, and  we can know that 
$\pi^{(i)}$ is an odd partition if $a$ is even and  $r\geq a$ or $r\leq a-1$ and $r$ is odd; $\pi^{(r)}$ is an even partition if $a$ is odd and $r\geq a$ or $r\leq a-1$ and $r$ is even.

We consider the lattice path that  corresponds to the partition $\lambda \in\overline{\mathcal{W}}_{N_1,\ldots,N_{k-1};a}$. By the bijection we constructed in Section 2, we know that:
\begin{itemize}
	\item For $r\geq a$ and $a$ is even,  the $n_r$ $r$-cluster correspond to $n_r$ peaks with the relative height $r$. To construct the $n_r$ peaks, we firstly add $n_r$ peaks at $(0,k-r-1)$ with weight $1,3,\ldots, 2n_r-1$ and relative height one.  
Then we add an SE step at $(0,k-r)$. Now  the $n_r$ peaks have even weight.  Under this condition, we konw that $\pi^{(r)}$ is a  partition with odd parts, therefore each peak increases in weight by an odd number.  During the construction, $r-a$ times  insertion of the SE step and $r-1$ times volcanic uplift will apply to these peaks, so  we get $n_r$ peaks whose weight have the same parity with $a$. 

\item	If $r\geq a$ and $a$ is odd. To construct the $n_r$ peaks  with relative height $r$, we firstly insert  $n_r$ peaks at $(0,k-r-1)$ with weight  $1,3,\ldots, 2n_r-1$ and relative height one.  
Then we add an SE step at $(0,k-r)$. Now  the $n_r$ peaks have even weight.  Under this condition, we know that $\pi^{(r)}$ is a partition with even parts, so each peak increases in weight by an even number.  During the construction, $r-a$ times  insertion of the SE step and $r-1$ times volcanic uplift will apply to these peaks, so  we get $n_r$ peaks whose weight have  the opposite parity with   $a$.  
 
 \item If $r\leq a-1$ and $r$ is odd, to construct the $n_r$ peaks with relative height $r$, we firstly add $n_r$ peaks at $(0,k-a)$ with weight $1,3,\ldots, 2n_r-1$ and relative height one.  
 Now  the $n_r$ peaks have odd weight. On this condition, $\pi^{(r)}$ is an odd partition, so each peak should increase in weight by an odd number. Now  the $n_r$ peaks have even weight. After  $r-1$ times volcanic uplift, we get $n_r$ peaks whose weight have the  opposite parity with the relative height $r$, then the weight are even.
 \item  If $r\leq a-1$ and $r$ is even then we construct the $n_r$ peaks with relative height $r$. Firstly, we add $n_r$ peaks at $(0,k-a)$ with weight  $1,3,\ldots, 2n_r-1$ and relative height one.  
 Now  the $n_r$ peaks have odd weight. In this condition, $\pi^{(r)}$ is an even partition. Then each peak increases in weight by an even number. Now  the $n_r$ peaks have odd weight.  After $r-1$ times volcanic uplift, we get $n_r$ peaks whose weight have the same parity with the relative height $r$, then the  weight of the peaks  are even.
\end{itemize}
Summerize the four conditions, we finished the proof of Theorem \ref{thmwp2}. \qed

\noindent{\bf The proof of Theorem \ref{thmle}}
 Aiming to give the proof, we proceed to construct the lattice paths which have  $n_r$ peaks with relative height $r$ and odd weight  satisfying the special $(k,a)$-conditions.	

For $r$ from $k-1$ to $1$, we construct the $n_r$ peaks successively:
\begin{itemize}
	\item [Step 1] Firstly, we insert  $n_r$ peaks with weight  $1,3,\ldots, 2n_r-1$ and relative height one at $(0,k-a)$, if $r\leq a-1$ or at $(0,k-r-1)$ otherwise. Then the major index of the lattice path increases $n_r^2+2n_r(n_{k-1}+\ldots+n_{r+1})$.

\item [Step 2]
If $r\geq a$, we add an SE step at $(0,k-r)$ which changes the parity of the peak's weight. In this step  the major index increases  $n_r+n_{r+1}+\cdots+n_{k-1}=N_r$. If $r\leq a-1$, we omit this step.

\item [Step 3] In this step, we consider four cases of $r$.
 
If $r\leq a-1$ and $r$ is even, we increase the weight of each peak with relative height one by one. Then the parity of the $n_r$ peaks become even as well as the other peaks do not change and  the major index of the lattice path increases $n_r$.

If $r\leq a-1$ and $r$ odd, we do nothing.

If $r\geq a$ and $a$ is odd, we increase the weight of each peak with relative height one by one. Then the parity of the $n_r$ peaks become odd as well as the other peaks do not change and  the major index of the lattice path increases $n_r$

If $r\geq a$ and  $a$ is even,
we do nothing.

\item [Step 4] Each peak with relative height one can increase in weight by an even number such that  each peak with relative height one increases a number no more than the peak with relative height one to its right. 

\item [Step 5] For $r$ from $k-1$ to $2$,  a volcanic uplift will be performed which changes the parity of the weight of the peaks and increase the major index by $1+3+5+\cdots+2N_r-1=N_r^2$.
\end{itemize}
Summerizing the above procedures and we can prove the peaks are all have even weight in either case as follows.
\begin{itemize}
\item [Case 1.] $r\leq a-1$ and $r$ even.
In the first step, the $n_r$ peaks with odd weight and then, after Step 3, the weight become even. In step 5, they will be performed $r-1$ times volcanic uplift which make them have the opposite parity with $r$, that is to say, the $n_r$ peaks have  odd weight.

\item[Case 2.]  $r\leq a-1$ and $r$ odd. By the  similar discussion of the first case,  one can see that the weight of the $n_r$ peaks  have the same parity with $r$, that is to say, the weight of these peaks are odd.

\item[Case 3.]  $r\geq a$ and $a$ odd. 
The first step get $n_r$ peaks with odd weight as well as the second step change the weight to be even, then increase weight by one makes the weight have odd parity.

The  $r-a$ times insertion of the SE step  and $r-1$ times volcanic uplift will change the partiy of weight $2r-a-1$ times, that is to say, an odd number changes parity $2r-a-1$ times, i.e. change even number of times becomes an odd number. Then the peaks have odd weight.

\item[Case 4.] $r\geq a$ and $a$ even, the similar discussion can show that in this case all peaks are also have odd weight. 
\end{itemize}
Then we have proved that our construction get a lattice path with odd weight peaks.

Now  we proceed to compute the generating function of the lattice path with $n_r$ peaks with relative height $r$.

For $k-1$ to $1$ the insertion of peaks increase the major index the following total number:
\[n_{k-1}^2+n_{k-2}^2+\cdots+n_1^2+2n_{k-2}N_{k-1}+\cdots+2n_1N_2;\]

For $r$ from $k-1$ to $2$, the volcanic lift increase the major index by the following number:
\[N_{k-1}^2+\cdots+N_2^2;\]
We combine the two forementioned numbers to get  
\[n_{k-1}^2+n_{k-2}^2+\cdots+n_1^2+2n_{k-2}N_{k-1}+\cdots+2n_1N_2+N_{k-1}^2+\cdots+N_2^2=N_1^2+\cdots+N_{k-1}^2;\]
For $r\geq a$ the insertions of the SE steps increase the major index by 
\[N_a+\cdots+N_{k-1}.\]

For $r\leq a-1$ and $r$ odd, operations in Step 3 increase major index by $n_r$, the total number increased  depend on the parity of $k,a$. 
 
In the case  $a$ is odd,  the major index increases:
\[n_2+n_4+\cdots+n_{a-1}+n_a+n_{a+1}+\cdots+n_{k-1}=n_2+n_4+\cdots+n_{a-3}+N_{a-1}.\] 
In the case  $a$ is even,  the major index increases: 
\[n_2+n_4+\cdots+n_{a-2}.\]

The increase in weight of each peak by any even number  are generated by the factor $1/(q^2;q^2)_{n_r}$. Then the generating function of the lattice path with $n_r$ peaks of relative height $r$ in the case   $a$ is odd is 
\[\frac{q^{N_1^2+N_2^2+\cdots+N_{k-1}^2+n_2+n_4+\cdots+n_{a-3}+N_{a-1}+N_{a}
		+N_{a+1}+\cdots+N_{k-1}}
}{(q^2;q^2)_{N_1-N_2}\cdots(q^2;q^2)_{N_{k-2}-N_{k-1}}(q^2;q^2)_{N_{k-1}}},\]
and  in the case   $a$ is even is 
\[\frac{q^{N_1^2+N_2^2+\cdots+N_{k-1}^2+n_2+n_4+\cdots+n_{a-3}	+N_a+N_{a+1}+\cdots+N_{k-1}}
}{(q^2;q^2)_{N_1-N_2}\cdots(q^2;q^2)_{N_{k-2}-N_{k-1}}(q^2;q^2)_{N_{k-1}}}.\]	
Then we have completed the proof of Theorem \ref{thmle}.

\noindent{\bf The proof of Theorem \ref{index}}	To prove this theorem we can show that given any  lattice path $L\in \mathcal{E}_{k,k}$,  each peak in $L$  has the same parity  with the  cluster corresponding to it in the parition $\lambda\in\mathcal{B}_{k,k}$ which  corresponds to $L$. We also can directly prove that the power of $y$ in \eqref{eqARRG} is the full even lower parity index of $L$. We employ the later one.

Now we recall the construction of a peak with relative height $r$ in a lattice path satisfying the special $(k,k)$-conditions.

All peaks with relative height $r$ were firstly inserted with odd weight and relative height one. At that time, the parity of peak    is odd. Because of  $a=k$, we do not need to insert any SE step at the start point. The insertion of other peaks and the $r-1$ times volcanic  uplift always change the parity of the weight and the relative height at the same time, that is to say the parity of all peaks are not change. and the lower even $r$-parity index and  full lower even parity index of the lattice path are not change.  Therefore the factor $q^{N_1^2+N_2^2+\cdots+N_{k-1}^2}$ generates the lattice paths in $\mathcal{E}_{k,k}$ with the lower even $r$-parity index and  full lower even parity index of the lattice path zero.

The factors in denominator generate even partitions which correspond to   even number increase in weight of each  peak in $L$ which do not change the parity of peaks in $L$. That is to say the lower even index of $L$ is still zero.

The factors $(-yq)_{n_r}$ generate  partitions with distinct parts no more than $n_r$, the power of $y$ denotes the number of parts in the partition. If $s$ is a part of the partition, we increase the weight of  rightmost $s$ peaks with relative height $r$ by one. Then the lower $r$-parity index increase one. So the number of parts in partition generated by $(-yq)_{n_r}$ is the lower $r$-parity index of lattice path. In conclusion, the full lower even  parity index is the power of $y$ in \eqref{eqARRG}.

By the discussion above, we have completed the proof of this theorem. \qed
	
\vspace{0.5cm}
	\noindent{\bf Acknowledgments.} This work was supported by  the National Natural Science Foundation of China.


\begin{thebibliography}{99} \small

\bibitem{agabre89}	A.K. Agarwal, D.M. Bressoud, Lattice paths and multiple basic hypergeometric series, Pacific J. Math. 135 (1989)  209–228.
		
		
\bibitem{and66} G.E. Andrews, An analytic proof of the Rogers--Ramanujan--Gordon
identities, Amer. J. Math. 88 (1966) 844--846.		
		
\bibitem{and74}G.E. Andrews, An analytic generalization of the Rogers--Ramanujan
identities for odd moduli, Proc. Nat. Acad. Sci. USA 71 (1974) 4082--4085.

		
\bibitem{and76}G.E. Andrews, The Theory of Partitions, Cambridge Mathematical Library, Cambridge University Press, Cambridge, 1998. Reprint of Addison-Wesley Publishing Co., 1976.
		
\bibitem{And10}
G. E. Andrews, Parity in partition identities, Ramanujan J. 23 (2010) 45-90.		

\bibitem{andbre84} G.E. Andrews, D.M. Bressoud, Identities in combinatorics III: Further aspects of ordered set sorting, Discrete
Math. 49 (1984) 222--236.

\bibitem{andbre85}	G.E. Andrews, D. Bressoud, On the Burge correspondence between partitions and binary words, in: Number Theory,
Winnipeg, MB 1983, Rocky Mountain J. Math. 15 (2) (1985) 225--233.
		
 		
\bibitem{bre79} D.M. Bressoud, A generalization of the Rogers--Ramanujan identities for all moduli, J. Combin. Theory Ser. A 27 (1979) 64--68.
		
\bibitem{Bre80} D.M. Bressoud, Analytic and combinatorial generalizations of the Rogers-Ramanujan identities, Mem. Amer. Math. Soc. 24 (1980) 1--54.
\bibitem{bre87}D. Bressoud, Lattice paths and the Rogers–Ramanujan identities, in: Number Theory, Madras 1987, in: Lecture Notes in Math., vol. 1395, Springer, Berlin, 1989, pp. 140--172.
		
		
		
\bibitem{chen13a}W.Y.C. Chen, D.D.M. Sang and D.Y.H. Shi, The Rogers--Ramanujan--Gordon Theorem for Overpartitions, Proc. London Math. Soc. 106 (3) (2013) 1371--1393.
		
\bibitem{chen13b}W.Y.C. Chen, D.D.M. Sang and D.Y.H. Shi, An Overpartition Analogue of Bressoud's Theorem of Rogers-Ramanujan Type, Ramanujan J., 36 (2015),  69–80, .
		
		
\bibitem{cor04}S. Corteel and J. Lovejoy, Overpartitions, Trans. Amer. Math. Soc. 356 (4) (2004) 1623--1635.

	
		
\bibitem{cor08}S. Corteel, J. Lovejoy and O. Mallet, An extension to overpartitions of the Rogers-Ramanujan identities for even moduli, J. Number Theory 128 (2008) 1602--1621 .
		
\bibitem{cor07}S. Corteel and O. Mallet, Overpartitions, lattice paths,
and Rogers–Ramanujan identities, J. Combin. Theory Ser. A 114 (2007) 1407–-1437.			
\bibitem{gor61}B. Gordon, A combinatorial generalization of the Rogers--Ramanujan identities,
		Amer. J. Math. 83 (1961) 393--399.
\bibitem{kim13}S. Kim and  A.J. Yee,  Rogers-Ramanujan-Gordon identities, generalized G\"{o}llnitz-Gordon identities, and parity questions, J. Comb. Theory, Ser. A 120(5), (2013) 1038--1056.
		
\bibitem{kur09}K. Kur\c{s}ung\"{o}z,
Parity considerations in Andrews--Gordon identities, European J. Combin. 31 (2010) 976--1000.

\bibitem{kur10}K. Kur\c{s}ung\"{o}z,
Cluster parity indices of partitions, Ramanujan J (2010) 23, 195–213.

\bibitem{lov03}J. Lovejoy, Gordon's theorem for overpartitions, J. Combin. Theory Ser. A 103 (2003) 393--401.



		
\end{thebibliography}
\end{document}